\documentclass[a4paper,11pt]{article}

\usepackage{setspace}
\usepackage[active]{srcltx}
\usepackage{makeidx}
\usepackage{multicol}
\usepackage{latexsym}
\usepackage{amsmath}
\usepackage{amssymb}
\usepackage{amscd}
\usepackage{theorem}
\usepackage[english]{babel}
\usepackage{xypic}

\usepackage[pdftex]{graphicx}
\usepackage{wrapfig}
\usepackage{capt-of}
\usepackage{pdfsync}

\newtheorem{theorem}{Theorem}[section]

\newtheorem{proposition}[theorem]{Proposition}
\newtheorem{lemma}[theorem]{Lemma}
\newtheorem{corollary}[theorem]{Corollary}
\newtheorem{definition}{Definition}[section]

\newtheorem{remark}[theorem]{Remark}

\newtheorem{remarks}[theorem]{Remarks}
\numberwithin{equation}{subsection}

\def\proof{\medskip \noindent \textit{Proof: }}
\def\qed{\hfill $\square$ }

\def\Or{O_{g_x}}
\def\d{\mathrm{d} }

\def\Ker{\mathrm{Ker}\,}
\def\VolX{\omega_X}

\def\proof{\medskip \noindent \textit{Proof:} }

\def\O{O_{g_x}}

\begin{document}

\title{On second-order, divergence-free tensors}

\date{}

\author{Jos\'{e} Navarro \thanks{Department of Mathematics, University of Extremadura,
Avda. Elvas s/n, 06071, Badajoz, Spain. \newline {\it Email address:} navarrogarmendia@unex.es \newline
The author has been partially supported by Junta de Extremadura
and FEDER funds.} }

\maketitle

\begin{abstract}

This paper deals with the problem of describing the vector spaces of divergence-free, natural tensors on a pseudo-Riemannian manifold that are second-order; i.e., that are defined using only second derivatives of the metric.

The main result establish isomorphisms between these spaces and certain spaces of tensors (at a point) that are invariant under the action of an orthogonal group.

This result is valid for tensors with an arbitrary number of indices and symmetries among them and, in certain cases, it allows to explicitly compute basis, using the theory of invariants of the orthogonal group.

In the particular case of tensors with two indices, we prove the Lovelock tensors are a basis for the vector space of second-order tensors that are divergence-free, thus refining the original Lovelock's statement.



\end{abstract}

\tableofcontents

\bigskip

\section*{Introduction}

Divergence-free tensors appear in a variety of places; among them, let us highlight that:

\medskip
- In a relativistic theory of gravitation, the field equation is assumed to be of the form:
$$ G_2(g) = T_2 $$  where $T_2$ is the energy-momentum tensor of matter, $g$ is the Lorentz metric of space-time that measures proper time and $G_2(g)$ is a suitable 2-tensor intrinsically constructed from the metric $g$; i.e., it is a {\it natural} 2-tensor.

By physical considerations (infinitesimal conservation laws), $T_2$ is assumed to be diver\-gence-free, so the field equation is determined by the choice of a di\-ver\-gen\-ce-free, natural tensor $G_2$. The paradigmatic example is the Einstein tensor, that fits into the field equation of General Relativity.

\medskip
- By the second Noether's Theorem, the  Euler-Lagrange tensor of a natural variational principle in the bundle of metrics is divergence-free. It is still an open problem to know whether any divergence-free, natural 2-tensor is the Euler-Lagrange tensor of a variational principle (\cite{AndersonPohj}, \cite{Takens}).

\medskip
In a pseudo-Riemannian manifold $(X,g)$, examples of natural tensors include tensor products and contractions among the metric, the dual metric, the curvature tensor or its covariant derivatives. The simplest natural tensors are those defined using only second derivatives of the metric, also called {\it second-order tensors}. 
This paper deals with the problem of describing those second-order tensors that are divergence-free. 

The major breakthrough in the area remains the work done by D. Lovelock regarding symmetric tensors with 2 indexes. In a couple of famous papers (\cite{LovelockI}, \cite{LovelockII}), he wrote down a basis for the $\mathbb{R}$-vector space of second-order 2-tensors that are symmetric and divergence-free. These Lovelock tensors are closely related to the existence of universal curvature identities (\cite{GilkeyI}) or the Chern-Gauss-Bonnet theorem (\cite{GilkeyLibro}), and have, since then, received much attention both by mathematical and physical communities (see, v. gr., \cite{Charmousis} or \cite{Kastor}). 

Nevertheless, a similar analysis for tensors with different symmetries or with a higher number of indices  revealed difficult (\cite{Anderson}, \cite{BedetLovelock}, \cite{Jaen}, \cite{Collinson}). Some progress has been made with the aid of computer programs (\cite{Jaen}), but there is still a lack of general results.

In this paper, we prove several statements, valid for tensors with any number of indices of symmetries, that allow to describe divergence-free tensors as tensors invariant under the action of an orthogonal group.

Firstly, we prove that divergence-free, second-order tensors are ``algebraic'', in the sense that their local expressions are polynomial functions of $g_{ab,cd}$: 

\medskip
{\bf Theorem \ref{PolinomialCharacter}:} {\it Any second-order tensor that is divergence-free is polynomial.}
\medskip

Its proof relies on some techniques introduced by Lovelock, that we develop so as to show the vanishing of derivatives of sufficiently large order. In the process, we make use of simple facts of graph theory, that allows to avoid lengthy calculations using symmetries of indices.

As a consequence, for example, it follows that any non-zero, second-order tensor that is divergence-free has an even number of indices.

\medskip

Secondly, it is easy to see that a polynomial tensor is divergence-free if and only if its homogeneous components are divergence-free (Proposition \ref{ReduccionHomog}). 

Then, we introduce certain spaces of tensors, $\mathrm{Div}^k$, depending on the symmetries and the number $p$ of indices under consideration, and we prove:

\medskip
{\bf Theorem \ref{TeoremaHomogeneos}} {\it Let $\Or$ be the group of linear automorphisms of $(T_xX , g_x)$, where $g_x$ is any pseudo-Riemannian metric with the fixed signature.

There exists an inclusion: 
$$ \left[ \begin{matrix} \mbox{  Second-order, natural $p$-tensors }   \\
\mbox{  homogeneous of degree $k$ and divergence-free  }
\end{matrix} \right] \,  \subset \, \left(  \mathrm{Div}^k \right)^{O_{g_x}}  \ . $$}
\medskip

That is to say, the problem is reduced to computing invariant tensors under the action of the orthogonal Lie group $O_{g_x}$, where the classical theory of invariants can be applied.

In the last Sections of the paper, we apply this technique to compute basis of divergence-free tensors in some particular cases.

In the simplest situation, that of tensors with 2 indices,  we recover Lovelock's result. Indeed, we prove a stronger statement, as we show that the Lovelock's tensors are also a basis for the vector space of second-order 2-tensors with null-divergence, but not necessarily symmetric:

\medskip
{\bf Theorem \ref{DivNulaOrdenDos}:} {\it The Lovelock tensors $ L_0 , \ldots , L_p $, where $2p \leq n-1$, are a basis for the $\mathbb{R}$-vector space of
second-order, natural 2-tensors that are divergence-free.}
\medskip

In dimension 4, this refined version was already established by Lovelock himself (\cite{LovelockII}), but the situation in higher dimensions remained open. 

On the other hand, as regards to tensors with some of their indices symmetric, we recover some of the results of and \cite{BedetLovelock} and \cite{Jaen}, and we prove a new statement for totally symmetric tensors:

\medskip
{\bf Theorem \ref{TotSymm}:} {\it Any second-order $2k$-tensor that is divergence-free is a constant multiple of the symmetrization of $g^* \otimes \stackrel{k}{\ldots} \otimes g^*$.}
\medskip

Finally, for differential forms we obtain a non-existence result:

\medskip
{\bf Theorem \ref{NoDifferentialForms}:} {\it There are no skew-symmetric, second-order tensors with zero-divergence, but the zero tensor.}
\medskip

\medskip
The paper is structured as follows:

The first preliminary Section collects together several classical definitions and results that are used later on.

Our later development heavily relies on the natural derivative of second-order tensors, originally developed by Lovelock. In Section 2, we present a new geometrical construction of this derivative, showing that naturalness corresponds to certain symmetries of the derivative.

Sections 3 and 4 are devoted to the proof of Theorems \ref{PolinomialCharacter} and \ref{TeoremaHomogeneos}, respectively. The last Sections contain the rest of the announced results.

\section{Preliminaries}

\subsection{Invariant theory of the orthogonal group}

Let $E$ be a finite dimensional $\mathbb{R}$-vector space and let $g$ be a non-singular metric of signature $(p,q)$ on it. Let $O_{p,q}$ 
be the orthogonal 
algebraic group of linear isometries of $(E,g)$. 

Later on, we make a strong use of the following classical theorem (\cite{Goodman}, \cite{Sancho}):  

\begin{theorem}\label{InvariantTheory} The vector space $\left( \otimes^r E^* \right)^{O_{p,q}} \, = \, \mathrm{Hom}_{O_{p,q}} ( \otimes^r E , \mathbb{R} ) $ of invariant linear forms on  $\otimes^r E$ vanishes if $r$ is odd, while, for $r$ even, it is spanned by contractions of the type:
$$ e_1 \otimes \ldots \otimes e_{r} \ \mapsto \ (g \otimes \ldots \otimes g) \, (e_{\sigma(1)} , \ldots  , e_{\sigma (r)} ) \   $$ where $\sigma$ is a permutation of $1, \ldots , r$.


\end{theorem}

Throughout the paper, this theorem is applied in conjunction with the semi-simplicity of the orthogonal group: if $\, V \subset \otimes^r E\,$ is a sub-representation of the group $\, O_{p,q}$, then any $O_{p,q}$-equivariant linear map $\, f \colon V \to \mathbb{R}$ is the restriction of a $O_{p,q}$-equivariant linear map $\, \widetilde{f} \colon \otimes^r E \to \mathbb{R}$. 

\begin{remark}
In the standard formulation of this result, the orthogonal group is considered as an {\it algebraic group}; that is, if the invariance of linear forms is understood as invariance under the action of all the points (real or imaginary) of the algebraic variety $O_{p,q}$.

In this paper the orthogonal group is considered as a {\it Lie group}; 
co\-rres\-pondingly, the invariance condition is understood as invariance under the action of the {\it rational (i.e., real) points} of the algebraic variety $O_{p,q}$ only.

Nevertheless, Theorem \ref{InvariantTheory} still holds in this setting (see \cite{Marco} or \cite{GilkeyII}). 
\end{remark}

\subsection{Second-order tensors}

Let $X$ be a smooth manifold of dimension $n$, and let $M \to X$ be the bundle of pseudo-Rie\-ma\-nnian metrics, with a fixed signature.

To fix notations, let us also recall that the {\it divergence} of a $p$-contravariant tensor $\,T\,$ on a pseudo-Riemannian manifold $\,(X, g)\,$ is the $(p-1)$-tensor:
$$ \mathrm{div}_g\, T := c_1^p (\, \nabla_g \, T ) \ ,  $$ 
where $c-1^p$ denotes the contraction of the first covariant with the $p^{th}$-contravariant indices. Its local expression, using summation over repeated indices, is:
$$ (\mathrm{div}_g\, T)^{j_1 \ldots j_{p-1}} = \nabla_k T^{j_1 \ldots j_{p-1} k} \ . $$

For simplicity, from now on we will only consider {\it contravariant} tensors:

\begin{definition}
A second-order $p$-tensor (not necessarily natural) is a second-order differential operator $\, T\colon M \rightsquigarrow \otimes^p TX\,$, that is, a morphism of bundles:
$$ T \colon J^2 M \to \otimes^p TX \  , $$ where $J^2M$ stands for the bundle of 2-jets of sections of $M$.
\end{definition}

In local coordinates, the components of a second-order $p$-tensor are smooth functions $$ T(g)^{j_1 \ldots j_p} \, = \, T^{j_1 \ldots j_p} (x_i , g_{ab}, g_{ab,c}, g_{ab,cd}) \ . $$ 

\smallskip
\begin{definition}  
The divergence of a second-order $p$-tensor $\,T \colon M \rightsquigarrow \otimes^p TX\,$ is the third order differential operator $\,\mathrm{div} \, T \colon M \rightsquigarrow \otimes^{p-1} TX\,$ defined as:
$$ (\mathrm{div}\, T) (g) := \mathrm{div}_g (T(g)) \ .  $$ 

A second-order $p$-tensor $T$ is divergence-free if $\mathrm{div}\, T$ is the zero map.
\end{definition}

\begin{proposition}\label{PropoLovelockLemma}
If a second-order $p$-tensor $T$ 
is divergence-free, then, on any chart, the functions $T^{j_1 \ldots j_p}$ satisfy:
\begin{equation}\label{SimetriaLinDivN} \frac{\partial T^{j_1 \ldots j_p }}{\partial g_{ab,cd}} +
 \frac{\partial T^{j_1 \ldots d}}{\partial g_{ab,j_p c}} + \frac{\partial T^{j_1 \ldots c}}{\partial g_{ab,d j_p}}  = 0 \ .
 \end{equation}
\end{proposition}

\proof In local coordinates:
\begin{align*}
(\mathrm{div}\, T)^{j_1 \ldots j_{p-1}} &= \frac{\partial T^{j_1 \ldots k}}{\partial x^k} 
+ \Gamma^{j_1}_{s k } T^{s \ldots k } + \ldots + \Gamma^{k}_{s k } T^{j_1 \ldots s} \\
& \, = \, \sum_{a\leq b} \sum_{c\leq d} \ \sum_{k=0}^n \frac{\partial T^{j_1 \ldots k}}{\partial g_{ab,cd}}\, g_{ab,cdk} +  F(x_i, g_{ab} , g_{ab,c}, g_{ab,cd} ) 
\end{align*}
for some smooth function $F$ on $J^2M$. 

Reordering this sum, we obtain:
\begin{align*}
0 \, = \, \sum_{a\leq b} \sum_{c\leq d \leq k} \   \left( \frac{\partial T^{j_1 \ldots k}}{\partial g_{ab,cd}} +
 \frac{\partial T^{j_1 \ldots d}}{\partial g_{ab,ck}} + \frac{\partial T^{j_1 \ldots c}}{\partial g_{ab,dk}} \right) g_{ab,cdk} + F(x_i , g_{ab} ,
 g_{ab,c} , g_{ab,cd})  \ ,
\end{align*} 
and the thesis follows because  $g_{ab,cdk}$ are elements of a chart on $J^3M$.

\qed

\subsubsection*{Second-order natural tensors}

The bundle $J^2M$ of jets of metrics is a {\it natural bundle}; that is, diffeomorphisms $\tau \colon U \to V$  between open sets of $X$ act on their sections:
$$ \tau \cdot j^2_x g := j^2_x ( \tau_* g) \ .   $$

\begin{definition}
A second-order tensor $T \colon J^2M \to \otimes^p TX$ is natural if it is a morphism of bundles that commutes with the actions of diffeomorphisms defined between open sets of $X$. 
\end{definition}

That is, $T$ is natural if, for any diffeomorphism $\tau \colon U \to V$ between open sets of $X$, it holds:
$$ T ( \tau \cdot j^2_x g ) = \tau_*  T(j^2_x g) \ . $$ 

\begin{definition}\label{definitionnormal}
The space of second-order normal tensors $N_2 \subset S^2 T^*_xX \otimes S^2 T^*_xX $ at a point $x$ is the kernel of the symmetrization $s_{3}$ in the last $3$ indices:
$$ 0 \to N_2 \to S^2 T^*_xX \otimes S^2 T^*_xX \xrightarrow{\ s_{3} \ } T^*_xX \otimes S^{3} T^*_xX \to 0 \ . $$
\end{definition}

\medskip
Any metric jet $j^2_x g$ defines a normal tensor $g^2_x \in N_2$: simply choose normal coordinates for $g$ at $x$ and write
$$ (g^2_x)_{ij kl} := g_{ij,kl} \ , $$ for the identity of the Gauss Lemma guarantees that the symmetrization of the last three indices of $g^2_x$ is zero.

\begin{lemma}\label{SymmetriesN} Second-order normal tensors have the following symmetries:
\begin{enumerate}
\item They are symmetric under the interchange of the first pair with the second pair of indices:
$$ G_{ij,kl} \ = \ G_{kl,ij} \ . $$

\item The cyclic sum of the last three indices is zero.
\end{enumerate}
\end{lemma}

\medskip
\noindent \textit{Proof:} Let us only check the first one, for the other is trivial. If $G \in N_2$, then:
\begin{align*}
G_{ij,kl} &= - G_{il,jk } - G_{ik,lj} = -G_{li,jk } - G_{ki,jl } \\
&= G_{lk,ij } + G_{lj,ki } + G_{kl,ij } + G_{kj,li } \\
&= 2 G_{kl,ij } + G_{jl,ki } + G_{kj,li } \\
&= 2 G_{kl,ij } - G_{ji,lk } - G_{jk,il } + G_{kj,li } = 2 G_{kl,ij } - G_{ij,kl } \ . 
\end{align*}

\qed

Local diffeomorphisms act transitively on $J^1M$; therefore, a natural tensor $T \colon J^2M \to \otimes^p TX$ is determined by its restriction to any fibre $\,\mathbb{A}_{j^1_x g}\,$ of the projection $J^2M \to J^1M$.

To compute the value of $T$ on an element $p \in \,\mathbb{A}_{j^1_x g}\,$ we can take normal coordinates; that is,  there exist a smooth map $\mathfrak{t}$ such that the following triangle commutes:
\begin{equation}
\xymatrix{
\mathbb{A}_{j^1_x g}\, \ar[r]^-{T} \ar[d]_-{\pi }  & \otimes^p T_xX  \\
N_2 \ar[ur]_-{\mathfrak{t}}   } \ 
\end{equation} where $\pi \colon A_{j^1_xg} \to N_2$ is the map that sends a 2-jet to its normal tensor.

These ideas lead to the following classical result (v. gr. \cite{Epstein}):

\begin{theorem}\label{Replacement} Let $g_x$ be a pseudo-Riemannian metric at $x$, and let $O_{g_x} $ be the orthogonal group of $(T_x X , g_x)$.

There exists an isomorphism of $\mathbb{R}$-vector spaces:
$$\begin{CD}
\left[ \phantom{\frac{1}{1}} \hskip-.2cm  \text{Natural  tensors } T \colon J^2M \to \otimes^p TX \right]  \\
@| \\
\left[ \phantom{\frac{1}{1}} \hskip-.2cm  \text{$\Or$-equivariant smooth maps } \mathfrak{t} \colon N_2\to \otimes^p T_xX \right]
 \end{CD} \ .  $$ 
\end{theorem}

Via this isomorphism, zeroth-order natural tensors corresponds with invariant elements in $\otimes^p T_xX $:
$$ \left[ \phantom{\frac{1}{1}} \hskip-.2cm  \text{Morphisms of natural bundles } T \colon M \to \otimes^p TX \right]  \, = \, \left( \otimes^p T_xX \right)^{\Or} \ . $$

\medskip

\section{Derivative of second-order tensors}

\medskip
The fibre $\,\mathbb{A}_{j^1_x g}\, $ of the projection $\,J^2_xM \to J^1_x M\,$ on $\,j^1_x g\,$ is an affine space, modelled on the vector space  $\,S_{2,2} := S^2 (T^*_xX) \otimes S^2 (T^*_xX)$. Hence, $\,S_{2,2}\,$ is the tangent space of $\,\mathbb{A}_{j^1_xg}\,$ at any point.

Let $\,\d \,$ denote the flat connection of the affine space $\,\mathbb{A}_{j^1_xg}$.

Let $\,T \colon J^2 M \to \otimes^p TX\,$ be a morphism of bundles. The restriction of $T$ to any of these fibres is a smooth map:
$$ T_{|} \colon  \mathbb{A}_{j^1_x g}  \ \to \ \otimes^p T_xX \ , $$ whose tangent linear map at a point $\,j^2_x g \in \mathbb{A}_{j^1_xg}\,$ is the tensor:
$$ T'_{j^1_xg} \, := \, \d _{j^2_xg} \left( T_{|} \right) \colon S_{2,2} \longrightarrow \otimes^p T_xX \ . $$

More generally, the $m^{th}$-covariant derivative defines the tensor:
$$ T^{m)}_{j^1_xg} \, := \, \d^m _{j^2_xg} \left( T_{|} \right) \colon S^m \left( S_{2,2} \right) \longrightarrow \otimes^p T_xX \ . $$


\begin{definition}
The derivative of $\,T \colon J^2M \longrightarrow \otimes^p TX\, $ is the morphism of bundles:
$$ T' \colon J^2 M \longrightarrow \otimes^p TX \otimes \left( S^2 TX \otimes S^2 TX \right) \quad , \quad j^2_x g  \, \longmapsto \, T'_{j^2_xg} \ .$$ 

Analogously, the higher derivatives are 
$$ T^{m)} \colon J^2 M \longrightarrow \otimes^p TX \otimes S^m \left(S^2 TX \otimes S^2 TX\right) \quad , \quad j^2_xg  \, \longmapsto \, T^{m)}_{j^2_xg} \ . $$ 
\end{definition}

In local coordinates, the coefficients of $T'$ are
\begin{equation}\label{ExpresionDerivada}
 T^{j_1 \ldots j_p; abcd}  \, = \, \frac{\partial T^{j_1 \ldots j_p}} {\partial g_{ab,cd} }  \ \ , 
 \end{equation} 
and, analogously, the coefficients of the $m^{th}$-derivative $\,T^{m)}\,$ are:
\begin{equation*}
T^{j_1 \ldots j_p; a_1b_1c_1d_1 \, \ldots \, a_mb_mc_md_m }  \, = \, \frac{\partial^m T^{j_1 \ldots j_p}} {\partial g_{a_1b_1,c_1d_1} \ldots \partial g_{a_mb_mc_md_m} }  \ . 
\end{equation*}

The local expression of the derivative
(\ref{ExpresionDerivada}) together with Proposition \ref{PropoLovelockLemma} imply:

\begin{corollary}[Lovelock]\label{LovelockLemma} If a second-order $p$-tensor $T$ 
is divergence-free, then its derivative $T'$ satisfies the following linear symmetry:
\begin{equation*}\label{LovelockId}
0 = \sum_{( j_p \, c \, d)} T^{j_1 \ldots j_p; abcd} := T^{j_1 \ldots j_p;abcd} +
T^{j_1 \ldots d;ab j_pc} + T^{j_1 \ldots c; ab dj_p} \ .
\end{equation*}
\end{corollary}

\subsubsection*{Natural tensors}


The naturalness of a second-order tensor is inherited by its derivative:

\begin{proposition}
The derivative of a second-order, natural  tensor is also a natural tensor.
\end{proposition}

\proof At any point $\,j^1_x g \in J^1M\,$, the naturalness of $\,T\,$ implies the commutativity of:
$$ 
\xymatrix{
\mathbb{A}_{j^1_x g} \ar[d]_{\tau_*} \ar[rr]^{\ T_{|} } &  &\otimes^p T_xX  \ar[d]^{\tau_* } \\
\mathbb{A}_{\tau_* j^1_x g}  \ar[rr]^{\ T_{|}}   & & \otimes^p T_{\tau(x)} X } \ . 
$$ 

Hence, the corresponding tangent linear maps at any point $\,j^2_xg \in \mathbb{A}_{j^1_xg}\,$ satisfy the commutative diagram:
$$ 
\xymatrix{ 
S_{2,2} \ar[d]_{\tau_*} \ar[rrr]^{\ T'_{j^2_xg} \ } & & &\otimes^p T_xX  \ar[d]^{\tau_* } \\
S_{2,2}  \ar[rrr]^{\ T'_{\tau_* (j^2_xg)} \ }  & & &  \otimes^p T_{\tau(x)}  X } \  
$$
that amounts to the naturalness of $\,T'\,$.  

\qed



\medskip
The space $\, N^2 \subset S^{2,2} := S^2 (T_xX) \otimes S^2 (T_xX)\, $ of contravariant normal tensors of order 2 is defined as the kernel of the symmetrization in the last 3 indices.

Let $\,j^1_x g \in J^1M\, $ and consider the map that sends a 2-jet to its normal tensor:
$$ \pi \colon \mathbb{A}_{j^1_x g} \longrightarrow  N_2 \qquad , \qquad j^2_xg \longmapsto g^2_x \  .  $$ 

This map is affine and the equations of its tangent linear map $\,\pi_*\,$ are:
\begin{equation}\label{AffineMap}
 S_{2,2} \xrightarrow{\ \ \pi_* \ \ } N_2 \quad , \quad G_{ij,kl} \ = \ \frac{1}{3} \left( 2S_{ij,kl}- S_{ik,jl} - S_{il,jk} \right) \ .  
 \end{equation}

Moreover, it is a  retract of the inclusion $\,N_2 \subset S_{2,2}\,$, and hence $\, S_{2,2} = N_2 \oplus \mathrm{Ker}\, (\pi_*) $.


\begin{lemma}\label{LemaIncidente} The subspace of $\,S^{2,2} = S_{2,2}^*\, $ incident with  $\,\mathrm{Ker}\, (\pi_*)\,$ is the space of contravariant normal tensors $\,N^2\,$.
\end{lemma}

\proof  Taking duals in $\, 0 \longrightarrow \mathrm{Ker}\, (\pi_*) \longrightarrow S_{2,2} \xrightarrow{\ \ \pi_* \ \ } S_{2,2}\, $, it follows that ${\rm Im}\, ( \pi^*)$ is the incident of $\Ker (\pi_*) $:
$$ 0 \longleftarrow (\Ker \pi_*)^* \longleftarrow S_{2,2}^* \xleftarrow{\ \ \ \pi^*  \ } S_{2,2}^* \ . $$ 

But $\, {\rm Im}\,  (\pi^*) = N^2\,$, because the dual map $\,\pi^* \colon S^{2,2} = S_{2,2}^* \to S_{2,2}^* = S^{2,2}\,$ can be checked to be defined by the same formula (\ref{AffineMap}), but for contravariant indexes.

\qed


\begin{proposition}[Symmetries of the derivative of natural tensors]\label{SymmetriesDer} If a  second-order tensor $\,T\colon J^2M \to \otimes^p TX\, $ is natural, then its derivative $\,T'\,$ takes its values in $\,\otimes^p TX \otimes N^2$:
\begin{equation} 
\xymatrix{
T' \colon  J^2M \ar[r]^-{} \ar[dr]_-{}  & \otimes^p TX \otimes S^{2,2}  \\
& \otimes^p TX \otimes N^2 \ar@{^{(}->}[u] }    \ .
\end{equation} 
\end{proposition}


\proof Let $\,j^2_xg \in J^2M\,$ and let $\,j^1_x g \,$ be its 1-jet, so that $\,j^2_x g \in \mathbb{A}_{j^1_x g}\,$. 

As $T$ is a natural tensor, there exists a $\Or$-equivariant smooth map $\,\mathfrak{t}\, \colon N_2 \longrightarrow \otimes^p T_xX\,$ such that 
$$ T\left( j^2_x \tilde{g} \right) \, = \, \mathfrak{t}\, \left(\tilde{g}^2_x \right) $$ for any $\,j^2_x \tilde{g}\,$ with the prefixed value $\,g_x\,$ at $\,x$.

In other words, the following triangle commutes:
\begin{equation}\label{Triangle}
\xymatrix{
\mathbb{A}_{j^1_x g} \ar[r]^-{T_{|}} \ar[d]_-{\pi }  & \otimes^p T_xX  \\
N_2 \ar[ur]_-{\mathfrak{t}}   } \ .
\end{equation} 

Hence, their tangent linear maps at $\,j^2_xg\,$ also commute:
\begin{equation*}
\xymatrix{
S_{2,2} \ar[rr]^-{T'_{j^2_xg} } \ar[d]_-{\pi_{*}}  & & \otimes^p T_xX  \\
N_2 \ar[urr]_-{\mathfrak{t}_*}   } \ 
\end{equation*}  proving that $T'_{j^2_xg } $ annihilates $\mathrm{Ker}\, \pi_* \subset S_{2,2} $. 

Therefore, via Lemma \ref{LemaIncidente} and the following isomorphisms $$\mathrm{Hom}_\mathbb{R} (S_{2,2} \, , \, \otimes^p T_xX) = \otimes^p T_xX \otimes (S_{2,2})^* = \otimes^p T_xX \otimes S^{2,2}$$ the map $\,T'_{j^2_xg}\, $ defines an element in $\,\otimes^p T_xX \otimes N^2$.

\qed

\section{Polynomial character of second-order, di\-ver\-gence-free tensors}

Let us begin this section with a short digression on graphs, that will be understood as finite $CW$-comple\-xes of dimension 1:

\begin{definition}
A graph is a compact Hausdorff topological space $K\,$, together with a finite subset $K_0 \subset K\,$, whose elements will be called {\it vertices}, such that:

\begin{enumerate}

  \item $K - K_0\,$ is a disjoint union of a finite collection of subspaces $e_i\,$, called {\it edges}, each of which is homeomorphic to an open interval. 
  
  \item The boundary of each edge is a pair of vertices or a single vertex.

\end{enumerate}
\end{definition}

Edges with equal endpoints will be called {\it loops}. Also, there can possibly be several edges between the same pair of vertices, v. gr:

\medskip
\centerline{\includegraphics[height=3cm]{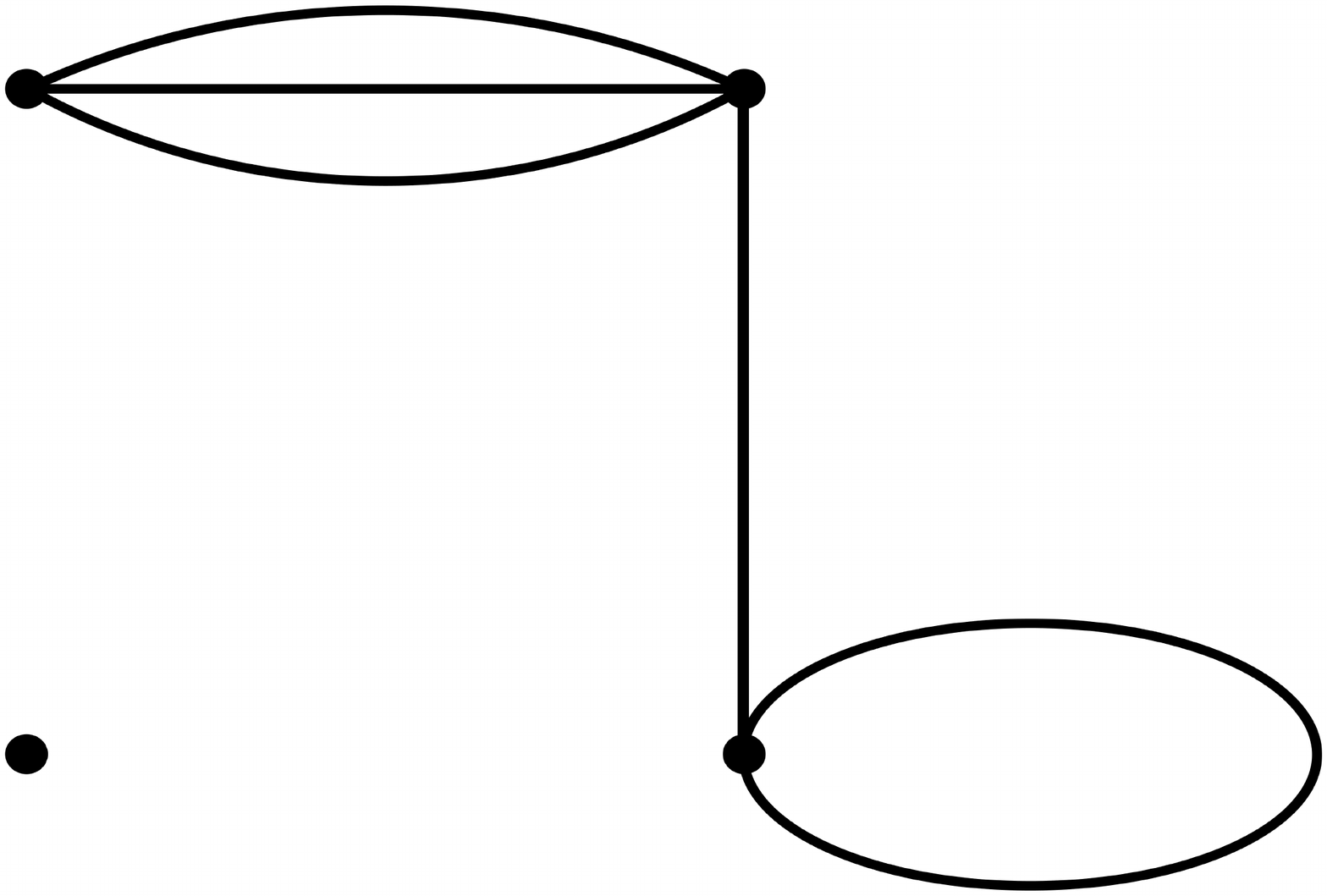}}
\medskip

Lattin letters $k,l,\ldots$ will represent vertices, and an edge joining the vertices $k$ and $l$ will be denoted as hyphened pair $k - l$.
The number of edges and vertices of a graph will be denoted $e$ and $v$, respectively .

\begin{definition} 
A cycle is a finite sequence of different edges $\,k_0 - k_1$, $\,k_1 - k_2, \ldots , k_m - k_0$, with $\,k_i\neq k_j\,$, for $\,i\neq j\,$.
A tree is a connected graph with no cycles.

\end{definition}

Any cycle satisfies the relation $\,e = v\,$, and any tree, the relation $\,e = v -1 \,$.


\begin{definition} A hair is a topological space homeomorphic to $\,[0,1)\,$, v. gr:

\medskip
\centerline{ \includegraphics[height=2.5cm]{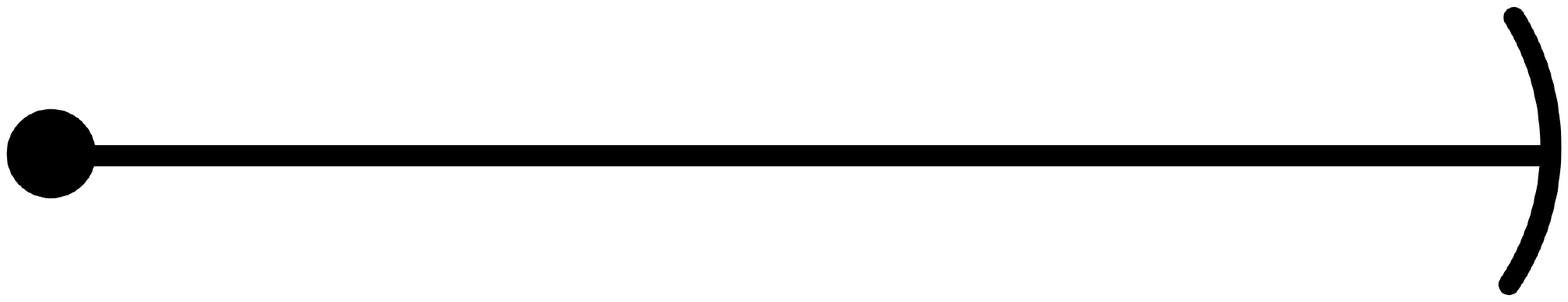} }
\medskip

A connected graph is a hairy cycle if it contains a cycle whose removal produces a disjoint union of hairs, v. gr:

\medskip
\centerline{ \includegraphics[height=3cm]{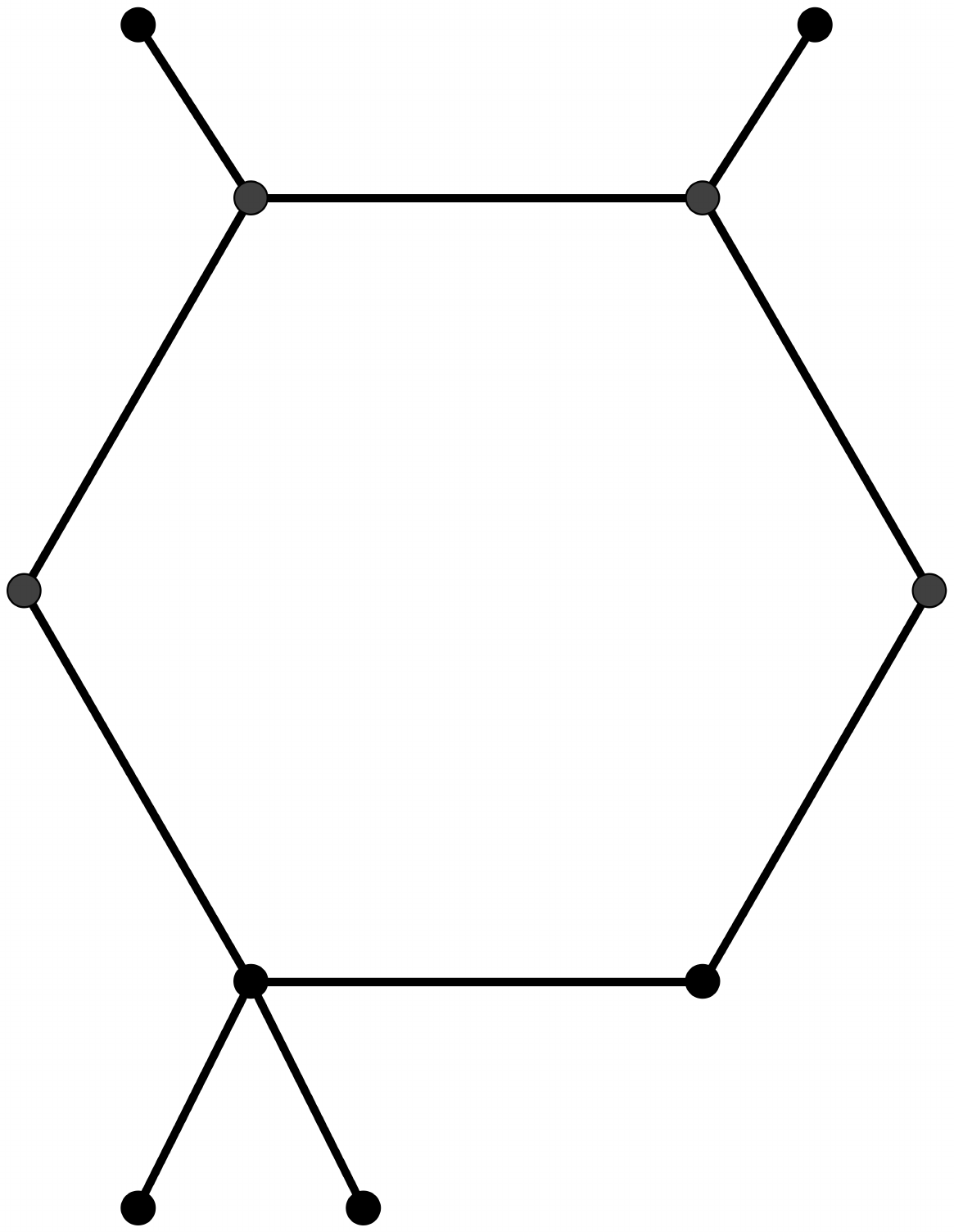} }
\medskip
\end{definition}

In particular, any cycle is a hairy cycle (the bald case, so to speak). 
Any hairy cycle also satisfies the relation $\,e = v\,$.

Finally, let us say a vertex is {\it simple} if there is only one edge arriving to it (the vertex of a loop is not considered to be simple), and that a vertex is {\it connected to a cycle} if there is a cycle in its connected component.

\begin{proposition}\label{KeyProposition} If a graph satisfies $\,e\geq v\,$, then one of the following options necessarily holds:
  \begin{enumerate}
      \item There exists an edge $\,k-l\,$ such that:
           \begin{itemize}
               \item[-] both $\,k\,$ and $\,l\,$ are not simple vertices;
               \item[-] after the removal of the edge $\,k-l\,$, the vertex $\,k\,$ is still connected with a cycle.
           \end{itemize}      
     \item The graph is a disjoint union of hairy cycles. 
\end{enumerate}
\end{proposition}

\proof If ${\it 1}$ is not satisfied, let us prove that the connected component of a cycle is a hairy cycle. 

Any vertex $\,h\,$ connected to the cycle, and not inside the cycle, has to be on an edge whose opposite endpoint is inside the cycle; otherwise there exists an edge $\,k-l\,$ satisfying ${\it 1}$ (see figure). 
Moreover, $\,h\,$ has to be simple; otherwise, the previous edge $\,k-h\,$  satisfies ${\it 1}$.
Finally, there does not exist an edge $\,k-l\,$ between two vertices of the cycle, different from the edges of the cycle, because such an edge would satisfy ${\it 1}$.

Therefore, the connected components of the graph are hairy cycles or trees.

\medskip
\centerline{ \includegraphics[height=3cm]{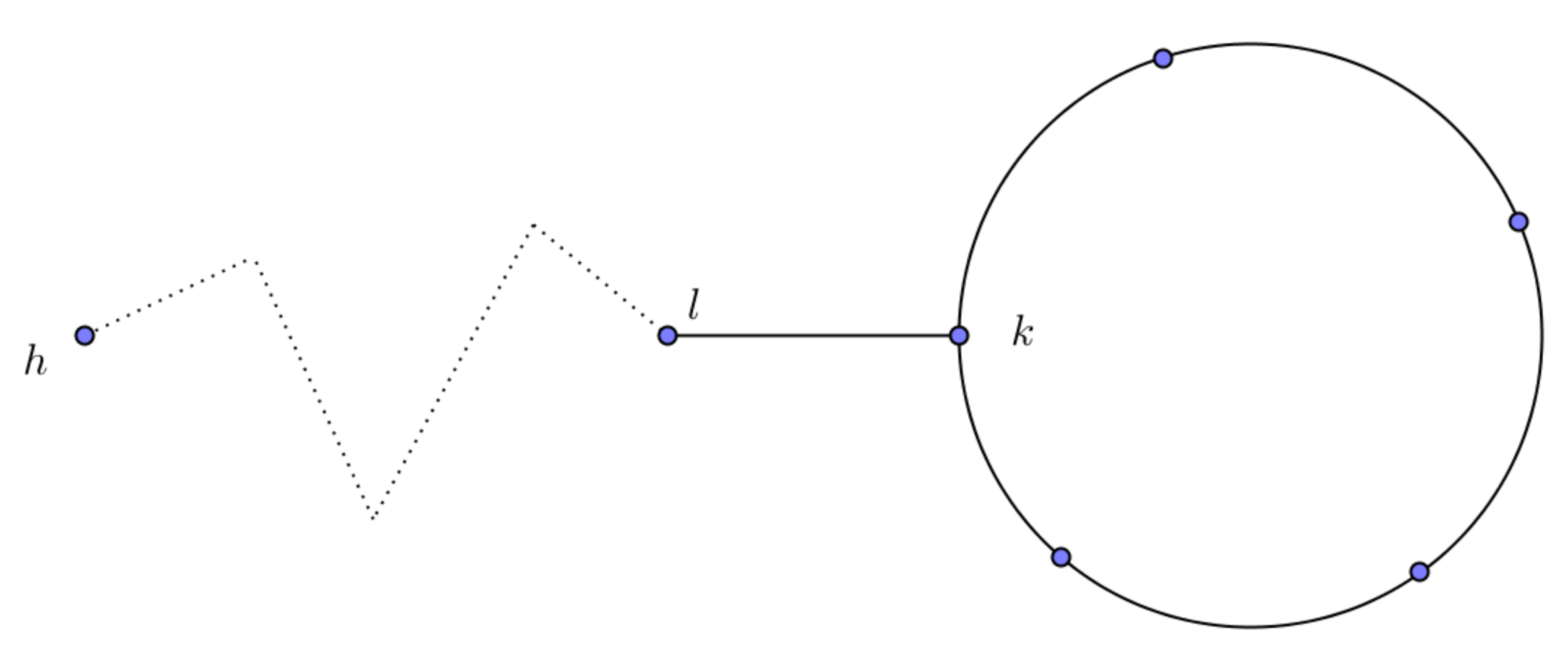} }
\medskip

As any hairy cycle satisfies $\,e = v\,$, and any tree $\,e = v -1\,$, the hypothesis $\,e \geq v\,$ for the graph implies no connected component is a tree.

\qed

\subsection*{Vanishing of derivatives}

\begin{definition} A second-order $p$-tensor $\,T \colon J^2M \to \otimes^p TX\,$ is polynomial (in the second deri\-va\-tives $\,g_{ab,cd}\,$) of degree $\,< m\,$ if its restrictions to the affine spaces $\,\mathbb{A}_{j^1_xg}$ $$ T_{|} \colon \mathbb{A}_{j^1_xg } \longrightarrow \otimes^p T_xX $$ are polynomial maps of degree $\,< m$, for any fibre $\,\mathbb{A}_{j^1_x g}$.
\end{definition}

This condition is equivalent to saying that the $m^{th}$-derivative $T^{m)}$ is null.

\medskip

Let $g_x$ be a pseudo-Riemannian metric at a point $x \in X$. By Theorem \ref{Replacement}, a second-order, polynomial natural tensor is defined by a polynomial $\O$-equivariant map $\mathfrak{t} \colon N_2\to \otimes^p T_xX$; that is, by an element in 

\begin{equation*}\label{Spaces}
 \bigoplus_{d_i} \mathrm{Hom}_{\O} ( S^{d_i} N_2 \ , \ \otimes^p T_xX ) \ . 
\end{equation*}

As these spaces of equivariant maps are spanned by total contractions of indices (Theorem \ref{InvariantTheory}) and tensors in $S^{d_i} N_2$ have an even number of indices, for any $d_i$, it readily follows:

\begin{proposition}\label{IndicesImpar}
If $p$ is odd, there are no polynomial, natural $p$-tensors $T \colon J^2 M \to \otimes^p TX$, but the zero tensor.
\end{proposition}


\begin{definition} For $m \geq 1$, let $\mathrm{Div}_x^m \subset \otimes^p T_xX  \otimes  S^m N^2  $
be the vector subspace whose elements satisfy:
\begin{align*}
 0 &= \sum_{(j_p \, c_1 \, d_1 )} T^{\, j_1 \, \ldots \,  j_p \, ; \, a_1 b_1 c_1 d_1 \,  \ldots } \ .
\end{align*}
\end{definition}

Due to the symmetries of $\,S^mN^2\,$ (any quatern $\,a_ib_ic_id_i\,$ can be put in the first position, and, by Lemma \ref{SymmetriesN}, $\,a_1b_1\,$ can be interchanged with $\,c_1d_1\,$), any element in $ \,\mathrm{Div}_x^m\,$ satisfies:
\begin{equation}\label{SimetriasDiv}
\sum_{ ( j_p a_i b_i ) } T^{\, j_1 \ldots j_p \, ; \, a_1 b_1 \, \ldots \, c_m d_m } \ = \  0 \  = \   \sum_{ ( j_p c_i d_i ) } T^{\, j_1 \ldots j_p \, ; \, a_1 b_1 \, \ldots \, c_m d_m }  \ . 
\end{equation} 

As a consequence, if the three indices $\,j_pa_i b_i\,$ are equal, then
$$ T^{\, j_1 \ldots l \, ; \,  \ldots \, ll \, \ldots } \ = \  0  \ . $$

\medskip
Due to Corollary \ref{LovelockLemma} and Theorem \ref{SymmetriesDer}, if a natural tensor $T \colon J^2M \to \otimes^p TX$ is diver\-gen\-ce-free, then its $m^{th}$- derivative takes value in this subspace:
\begin{equation}\label{Factorizacion}
\xymatrix{ T^{m)}\colon J^2_x M \ar[r] \ar[dr] & \otimes^p T_xX \otimes S^m N^2 \\
& \mathrm{Div}_x^m \ar@{^{(}->}[u] } \  . 
\end{equation}

Therefore, 

\begin{proposition} 
If there exist $\,m \in \mathbb{N}\,$ such that $\,\mathrm{Div}_x^m = 0$, then any diver\-gen\-ce-free, natural $p$-tensor $\,T \colon J^2M \longrightarrow \otimes^p TX\,$ has to be polynomial (in the second derivatives of the metric), of degree less than $m$.
\end{proposition}





Consider a component of an element in $\mathrm{Div}_x^m$:  
$$ \mathsf{T} = T^{\, j_1 \ldots j_p \, ; \, a_1 b_1 \, \ldots \, c_md_m  } \  , $$ that we understand as a linear function on $\mathrm{Div}_x^m$.

\begin{definition}
Its associated graph is defined as follows:
\begin{itemize}
\item It has $\,n\,$ vertices, labelled by $\,k \in \{ 1, \ldots , n \}$, where $\,n = \dim X$.

\item For each pair of indexes $\,a_ib_i\,$ (or $\,c_id_i\,$) in $\,T^{\, j_1 \ldots j_p \, ; \, a_1 b_1 \, \ldots \, c_md_m  }$, there is one edge $\,a_i-b_i\,$ joining the vertices $\,a_i\,$ and $\,b_i\,$ (resp. an edge $\,c_i - d_i\,$ joining $\,c_i\,$ and $\,d_i\,$).
\end{itemize}
\end{definition}

The indexes $\,j_1 \ldots j_p\,$ are irrelevant to construct the graph. As an illustration, consider the following example, in which  $\,n=6 \,$ and $\,m=3\,$:
\begin{equation}\label{Example}
\mathsf{T} = T^{\, j_1 \ldots j_p \, ; \, 11  43 \, ; \,  12  12 \, ; \, 62  61} \ , 
\end{equation} 
whose associated graph is: 


\medskip
\centerline{\includegraphics[height=3cm]{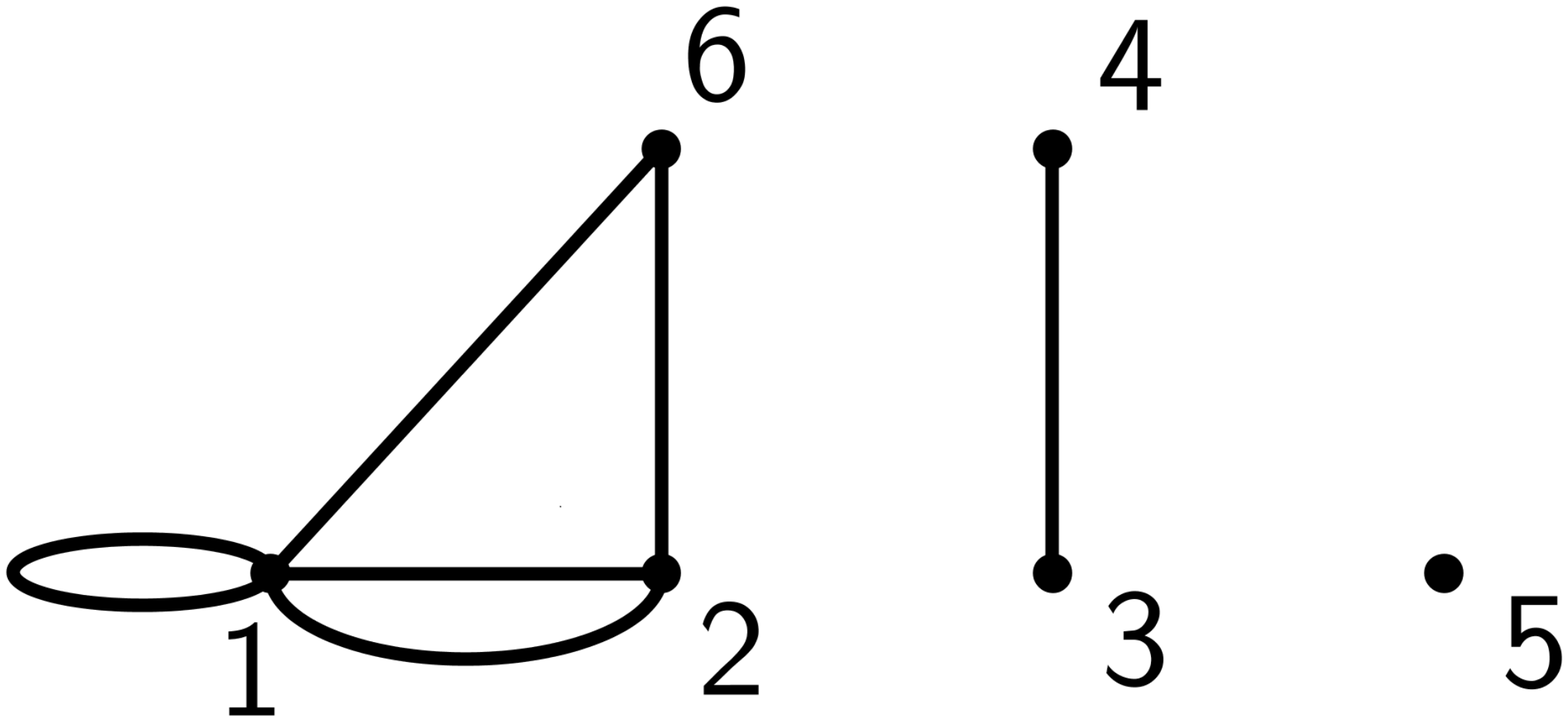}}
\medskip


\begin{lemma}\label{LemmaLoop}
Consider  a component of an element in $\mathrm{Div}_x^m$:
$$ \mathsf{T} = T^{\, j_1 \ldots j_p \, ; \, a_1 b_1 \, \ldots \, c_md_m  } \  . $$ 

If, in the associated graph, the vertex $\,j_p \in \{ 1 , \ldots , n \}\,$ is connected to a cycle, then this component is zero.
\end{lemma}

\proof By hypothesis, there exist edges $j_p-k_1$, $k_1-k_2$, $\ldots$, $k_r - l$ connecting $j_p$ with a cycle $l- m_1$, $\ldots$, $m_s - l$. 

Iterated use of (\ref{SimetriasDiv}) yields ($\sim$ denotes proportional by a non-zero factor):
\begin{align*}
T^{\, \ldots j_p \, ; \, \ldots \, j_p k_1 \ldots   } & \, \sim \, T^{\,  \ldots k_1 \, ; \, \ldots \, j_p j_p \ldots k _1 k_2 \ldots   } \, \sim \, T^{\,  \ldots k_2  \, ; \, \ldots k_1k_1 \ldots k_2k_3 \ldots   } \,  \sim \ldots  \\
& \, \sim \, T^{\,  \ldots l \, ; \, \ldots \, l m_1 \ldots   } \sim \, T^{\,  \ldots m_1 \, ; \, \ldots l l \ldots m_1 m_2 \ldots   }  \sim \, T^{\, \ldots m_2 \, ; \, \ldots \, l l  \ldots m_2 m_3 \ldots  }  \,  \sim \,
\ldots \,  \\
&  \sim \, T^{\, \ldots l \, ; \, \ldots \, l l  \ldots   } = 0 \ .
\end{align*} 
\qed

\begin{theorem}\label{PolinomialCharacter}  If $\,m \geq n/2$, then: 
$$\mathrm{Div}_x^m = 0 \ . $$

As a consequence, any divergence-free, natural tensor $T \colon J^2M \to \otimes^p TX$ is polynomial (in the second derivatives of the metric), of degree less or equal than $(n-1)/2$.
\end{theorem}

\proof  Let $\,m \geq n/2\,$ and suppose there exist a component $\,\mathsf{T}\, $ which does not vanish identically on $\,\mathrm{Div}_x^m$.

Another component $\,\mathsf{T}'\,$ is considered equivalent to $\,\mathsf{T}\,$ if there exists $\,\lambda \neq 0$ such that $ \,\mathsf{T} = \lambda \mathsf{T}'\,$ (as linear functions on $\,\mathrm{Div}_x^m$) and if the last $4m +1$ indices of $\,\mathsf{T}'\,$ are a permutation of those in $\,\mathsf{T}$.

Among all the elements  equivalent to $\,\mathsf{T}\,$, let 
$$ \mathsf{T}_R = T^{\, j_1 \ldots \,  j_p \, ; \, a_1 b_1  \ldots c_m d_m }$$  be one with the highest number of loops (i.e., edges $\,k-k\,$ with equal endpoints) in the associated graph. As $2m \geq n$, the graph associated to $\,\mathsf{T}_R\,$ satisfies $ \,e \geq v\,$, and we can invoke Proposition \ref{KeyProposition}. 

If the graph associated to $\,\mathsf{T}_R\,$ is a disjoint union of hairy cycles, then any vertex is connected with a cycle, and hence $\,\mathsf{T}_R = 0\,$ (Lemma \ref{LemmaLoop}), in contradiction with the hypothesis.

Otherwise, let $\,k-l\,$ be an edge as in {\it 1} in Proposition \ref{KeyProposition}. Then:
$$ T^{\,  \ldots j_p \, ; \, \ldots \, k l  \ldots } \, = \, - T^{\,  \ldots k \, ; \, \ldots \, j_p l \ldots   } -  T^{\,  \ldots l \, ; \, \ldots \, k j_p \ldots   } \ . $$

The first addend is zero because $\,k\,$ is connected to a cycle in the graph of  $\,T^{\,  \ldots k \, ; \, \ldots \, j_p l ...}$. 

If $\,l=k\,$,  the second addend is also zero and thus $\,\mathsf{T}_R = 0$.

In other case,  $\,l \neq k\,$ and, as $\,l\,$ is not simple, there exists at least another edge $\,l -m\,$ in the graph of $\,T^{\,  \ldots l \, ; \, \ldots \,  k j_p \ldots   } $. If $\,m = l$, then $\,\mathsf{T}_R = 0$, because:
$$ T^{\,  \ldots l \, ; \, \ldots \, l l  \ldots   } \, = \, 0 \ , $$
and, if $\,m \neq l$, we also arrive to a contradiction, because we produce a component equivalent to $\,\mathsf{T}_R\,$ but with a greater number of loops:
$$  T^{\,  \ldots j_p \, ; \, \ldots \, k l  \ldots } \, \sim \, T^{\, \ldots l \, ; \, \ldots \, l m  \ldots   } \, \,  \sim \, \,  T^{\,  \ldots m \, ; \, \ldots \, l l  \ldots   } \ .  $$
\qed

\section{Computation of divergence-free tensors}

Let us explain how the results in the precedent Section allow to reduce the computation of second-order $p$-tensors that are divergence-free to a problem of invariants for the orthogonal group.
\medskip

Let $\,T\colon J^2M\longrightarrow\otimes^pTX\,$ be a second-order, natural $p$-tensor that is di\-ver\-gence-free. By Theorem \ref{PolinomialCharacter}, this tensor is polynomial on the second derivatives of the metric, of degree $m\leq (n-1)/2$. Let us decompose it as
$$T\,=\, T_0+\cdots+T_m\ ,$$
where the  tensors $T_k$ are homogeneous polynomials of degree $k$ in the second de\-ri\-vatives of the metric. As any diffeomorphism of the base manifold $X$ acts linearly on the coordinates $g_{ab,cd}$ of $J^2M$, it is trivial to check that each tensor $T_k$ is natural. 

Let $\,g_x\,$ be a fixed metric at a point $\,x\in X$.  
By Theorem \ref{Replacement}, the natural tensor $\,T\colon J^2M \to \otimes^pTX$ corresponds with a polynomial $\Or$-equivariant map $\,\mathfrak{t}\colon N_2 \to  \otimes^pT_xX $, satisfying $\, T(j^2_xg)\, = \, \mathfrak{t}\, (g^2_x)\,$ on each metric jet $j^2_x g$ having the prefixed value at $x$.

 If $\,\mathfrak{t}\,= \,\mathfrak{t}_0 \,+\ldots +\, \mathfrak{t}_m$ is the decomposition into homogeneous components, then it is easy to check that each addend $\,\mathfrak{t}_k\,$ corresponds to the $p$-tensor $\,T_k$. This shows that the natural $p$-tensor $\,T_k\,$ is homogeneous of weight $\,-2-2k\,$ respect to the metric (i.e., it satisfies $T_k (\lambda^2 g) = \lambda^{-2-2k} T_k (g)$) and, consequently, so does its divergence $\,\text{div}\, T_k$. 

Therefore, in the equality 
$$0\,=\,\text{div}\, T\,=\, \text{div}\, T_0+\cdots+\text{div}\, T_m$$
each addend has a different weight, so each $\,\text{div}\, T_k\,$ is zero.

All in all, we have proved:
 
\begin{proposition}\label{ReduccionHomog}
Any natural $p$-tensor $T \colon J^2 M \to \otimes^p TX $ that is divergence-free admits a decomposition
$$T\,=\, T_0+\cdots+T_m\ ,$$
where each addend $T_k$ is a divergence free, second-order, natural $p$-tensor which is a homogeneous polynomial of degree $k$ in the second derivatives of the metric.
\end{proposition}

\medskip
As regards to the computation of homogeneous tensors, consider the 
commutative triangle: 
\begin{equation*}
\xymatrix{
\mathbb{A}_{j^1_x g} \ar[rr]^-{T_{|}} \ar[d]_-{\pi } &  & \otimes^p T_xX  \\
N_2 \ar[urr]_-{\mathfrak{t}}   } \ ,
\end{equation*} where $\,\pi(j^2_xg)=g^2_x$, and $j^1_x g$ is any jet with the prefixed value at $x$. 

As  $\,\pi\,$ is a surjective affine map, $\,T_|$ is a homogeneous polynomial of degree $\,m\,$ if and only if so it is $\, \mathfrak{t}$. Therefore, the following bijections hold:

$$\begin{CD}
\left\{
\begin{aligned}
&\text{Second-order natural tensors } \ \  T \colon J^2M \longrightarrow\otimes^pTX \\
& \, \text{homogeneous of degree } m \text{ in the second derivatives} \, 
\end{aligned}\right\}\\
@| \\
\left\{
\begin{aligned}
& \Or \text{-equivariant polynomials }\ \mathfrak{t}\colon N_2\longrightarrow\otimes^pT_xX \ \\
&\hskip 1.8cm \text{ homogeneous of degree } m
\end{aligned} \right\}\\
@| \\
\text{Hom}_{O_{g_x}} \left(S^mN_2\, ,\,\otimes^pT_xX \right)\,=\, \left( \otimes^pT_xX\,\otimes\, S^mN^2\right)^{O_{g_x}}\\
@| \\
\left\{ 
\begin{aligned}
&\text{Zeroth-order natural  tensors } \\
& \hskip .3cm \tilde T\colon M\longrightarrow\otimes^pTX\otimes S^mN^2 
\end{aligned}
\right\}
\end{CD}$$

We use the symbol $\,N^2\,$ to denote the space of contravariant normal tensors at $\,x$, as well as to denote the bundle of such tensors.

This sequence of maps sends each natural tensor $\, T\colon J^2M\to\otimes^pTX$, homogeneous of degree $m$ in the second derivatives of $\,g$, to its $\,m^{th}$-derivative $\,T^{m)}\colon M\to \otimes^pTX\otimes S^mN^2$. 

Observe that, after differentiating $\,m\,$ times, the tensor $\,T^{m)}\,$ does no longer depend on the second derivatives of the metric; hence, by naturalness, nor does it on the first derivatives.


The divergence-free condition can now be imposed, resulting:

\begin{theorem}\label{TeoremaHomogeneos} Let $\,g_x\,$ be a metric at a point $\, x\in X$. For any $\,m \geq 1\, $ there exists an injective map
$$\begin{CD}
\left\{
\begin{aligned}
&\text{Divergence free, natural tensors }\ T \colon  J^2M\longrightarrow\otimes^pTX\\
&\text{homogeneous  of degree }m\text{ in the second derivatives}
\end{aligned}\right\}\\  |\bigcap \\  
 \quad\left(\mathrm{Div}_x^m\right)^{O_{g_x}}
\end{CD}$$ that sends a tensor $T$ to its $m^{th}$-derivative $\,T^{m)}\,$  at the point $x$.
\end{theorem}

\proof If the tensor $\,T\,$ is divergence-free, then its $m^{th}$-derivative (at $x$),   $\, T^{m)}\,$ takes its values inside the subspace  $\,\text{Div}_x^m \subset\otimes^pT_xX\otimes S^mN^2$, by formula (\ref{Factorizacion}).

And the vector space of zeroth-order natural tensors $\widetilde{T} \colon M \to \otimes^p TX \otimes S^m N^2$ that take values inside $\mathrm{Div}_x^m$ is isomorphic to $\left(\mathrm{Div}_x^m\right)^{O_{g_x}}$, by Theorem \ref{Replacement}.

\qed


\begin{remark}\label{RemarkEquivalencia}
Although it will not be needed in the rest of the paper, let us remark that this inclusion is indeed an isomorphism.

This is due to the fact that a second-order, natural tensor $T$ is divergence-free {\it if and only if} its local components satisfy (\ref{SimetriaLinDivN}). This statement can, in turn, be proved using the formula:
$$ \nabla_k T^{j_1 \dots j_p} \ = \ \frac{2}{3} \sum_{a,b,c,d} \frac{\partial T^{j_1 \dots j_p}}{\partial g_{ab,cd}} \nabla_k R_{cabd} $$ which
can be found, for example, in \cite{DuPlessis}.
\end{remark}

\section{Lovelock tensors}

In this Section, let us consider tensors with $p=2$ indices. In this case, 
$$ \mathrm{Div}_x^m  \subset \otimes^2 T_x X \otimes S^m T^*_x X $$ is the vector subspace defined by the equations:
\begin{align*}
 0 &= \sum_{(j \, c_1 \, d_1 )} T^{\, i j  \, ; \, a_1 b_1 c_1 d_1 \,  \ldots } \ .
\end{align*}

Let $g_x$ be a pseudo-Riemannian metric at a point $x \in X$. The key computation is the following:

\begin{lemma}\label{CalculoDimension}
For any $m \geq 1$: 
$$ \dim \left( \mathrm{Div}_x^m \right)^{\Or} \leq 1 \ . $$ 
\end{lemma}

\proof The vector space $\, \left( \mathrm{Div}_x^m \right)^{\Or} \, = \, \mathrm{Hom}_{\Or} \left( \mathbb{R} \, , \, \mathrm{Div}_x^m \right)\,$ is isomorphic to:
$$ \mathrm{Hom}_{\Or} \left( \left( \mathrm{Div}_x^m\right)^* \, , \, \mathbb{R} \right) $$ which, in turn, is spanned by iterated contraction of indices (Theorem \ref{InvariantTheory}).

Let us prove that any total contraction of indices is proportional to:
$$  T^{ii; jj  \ldots kk } \ , $$ where equal letters denote contraction of the corresponding positions.

We argue by descendent induction on the number of contracted pairs (i.e., contraction of an index in an odd position with the index in the following position).

Given a total contraction, if the first index is not contracted with the second one, then:
$$ T^{ij \, ; \, \ldots \, j k \ldots } \, \sim \, T^{i k  \, ; \, \ldots \, jj \ldots } $$ and the induction hypothesis applies.

Otherwise, we can assume the third index is not contracted with the fourth, and hence:
$$ T^{ii \, ; \, jk \ldots  } \, = \, - T^{ij \, ; \, ik\, \ldots \, j m  \ldots }  -  T^{i k \, ; \, i j \, \ldots \,  k l  \ldots }
 \, \sim  \,  T^{im \, ; \, \ldots \, jj  \ldots } + T^{il \, ; \, \ldots \, kk  \ldots } \ . 
 $$ This two addends have the same number of contracted pairs as the original one, so we are reduced to the previous case. 
 
 \qed

\medskip
As a consequence, the vector space of divergence-free 2-tensors, homogeneous of degree $k \leq (n-1)/2$, has dimension at most one. In the following subsection, let us define explicit generators for these spaces.

\subsubsection*{Definition of the Lovelock tensors}

Let $g$ be a pseudo-Riemannian metric and let us consider it as a one-form with values on one-forms.


Its Riemann-Christoffel tensor $R$ can also be understood
as a 2-form with values on 2-forms that is symmetric, i.e., a section of $S^2 \left( \Lambda^2 T^* X \right) \subset \Lambda^2 T^*X \otimes \Lambda^2 T^*X$.


With this language, the differential Bianchi identity and the torsion-free property of the Levi-Civita connection $\nabla$ amount to the equations:
\begin{equation}\label{BianchiIdentity} \d_\nabla R = 0 \quad , \quad \d_\nabla g = 0  \ .  \end{equation}

With respect to the wedge product of forms, consider the following $(n-1)$-forms with values on $(n-1)$-forms:
\begin{align*}
\widetilde{L}_{k} &:= R \, \wedge \, \stackrel{k}{\ldots} \, \wedge \, R \, \wedge g \, \wedge \stackrel{n- 2k -1}{\ldots} \wedge \, g
\end{align*} where $k$ runs from 0 to the integer part of $(n-1)/2$.

These $\widetilde{L}_k$ are clearly symmetric, i.e., sections of $S^2 (\Lambda^{n-1} T^*X) \subset \Lambda^{n-1} (T^*X) \otimes \Lambda^{n-1}(T^*X) $, and also satisfy $\d_\nabla \widetilde{L}_k = 0$, in virtue of (\ref{BianchiIdentity}).

The following statement is easy to prove (see, v. gr., \cite{LovelockAl}):

\begin{proposition}\label{LovelockDuals}
Contraction with a volume form, $\VolX $, defines a linear
isomorphism:
$$ TX \otimes TX \ \xrightarrow{\sim} \ \Lambda^{n-1} (T^*X) \otimes \Lambda^{n-1}(T^*X) \quad , \quad D \otimes D' \ \mapsto \ i_{D} \VolX \otimes i_{D'} \VolX \ , $$
and symmetric 2-tensors correspond with sections of $S^2 (\Lambda^{n-1} T^*X)$.

Moreover, if $T$ and $\Pi$ are a 2-tensor and a valued $(n-1)$-form corresponding via this isomorphism, then:
$$ \d_{\nabla} \, \Pi = 0 \quad \Leftrightarrow \quad \mathrm{div}\, T = 0 \ . $$
\end{proposition}

\begin{definition}\label{DefinitionLovelock} The Lovelock's tensors $\,L_k\,$ are the 2-contravariant tensors on $\,X\,$ corresponding to the vector-valued forms $\,\widetilde{L}_k\,$ via the isomorphism above.

Hence, they are symmetric and divergence-free 2-tensors.
\end{definition} 

\medskip
Apart from the trivial case of the dual metric, $L_0 = g^*$, the simplest Lovelock tensor, $L_1$, is proportional to the contravariant Einstein tensor; i.e., via the isomorphism above,
\begin{equation*}
 R\wedge g\wedge\overset{n-3}{\dots}\wedge g  \quad  \longmapsto  \quad  (-1)^{q+1} \,(n-3)!\,\left( Ric -\frac{r}{2}g^* \right)
\end{equation*}
where $q$ stands for the number of $-1$ in the signature $(p,q)$ of $g$.

In general, it is easy to check that the $k^{th}$-Lovelock tensor $L_k$ is a homogeneous polynomial of degree $k$ on the second-derivatives of the metric. Therefore, $\,L_k\,$ generates the vector space of divergence-free tensors $T \colon J^2M \to \otimes^2 TX$  that are homogeneous polynomials of degree $k$ on the second derivatives of the metric.

To be precise:

\begin{theorem}\label{DivNulaOrdenDos} The Lovelock tensors $ \,L_0 , \ldots , L_m $, where $\,2m \leq \dim X -1$, are a basis for the $\mathbb{R}$-vector space of
second-order, natural 2-tensors that are divergence-free.

That is,
$$\left[ \begin{matrix} \mbox{ Natural tensors  } T \colon J^2M \to \otimes^2 TX  \ \\
\mbox{ that are divergence-free }
\end{matrix} \right] \ = \ \langle L_0 , \ldots , L_m \rangle \ . $$
\end{theorem}

\proof
By Theorem \ref{PolinomialCharacter}, any tensor $\,T\,$ of the type under consideration is polynomial, of degree $\,m\leq (n-1)/2\,$ in the second derivatives of the metric. 

By Proposition \ref{ReduccionHomog},  $\,T=T_0+\cdots+T_m$, where each $\,T_k\,$ is divergence-free and is homogeneous of degree $k$ in the second derivatives of the metric. 

As the space of such tensors has dimension  $\,\leq 1$, due to Theorem \ref{TeoremaHomogeneos} and Lemma \ref{CalculoDimension}, the tensor $\,T_k\,$ coincides, up to a constant factor, with $\,L_k$. 

\qed

\begin{remarks} 
Our proof also characterizes the $k^{th}$-Lovelock tensor $\,L_k\,$ as the only, up to a constant factor, second-order, natural $2$-contravariant tensor which is divergence-free and homogeneous of weight $w=-2-2k$.
\end{remarks}

\section{Other computations}

Let $p \geq 4$ (recall that, due to Theorem \ref{PolinomialCharacter} and Proposition \ref{IndicesImpar}, there are no second-order, divergence-free natural tensors with an odd number of indices, but the zero tensor).

As an illustration of the techniques explained above, let us firstly consider tensors that are symmetric in 3 indices; i.e., sections of:
$$ \otimes^{p-3} T_xX \otimes S^3 T_xX \ . $$

\begin{definition} For $m \geq 1$, let $S\mathrm{Div}_x^m \subset \otimes^{p-3} T_xX \otimes S^3 T_xX  \otimes  S^m N^2  $
be the vector subspace whose elements satisfy:
\begin{align*}
 0 &= \sum_{(j_3 \, c_1 \, d_1 )} T^{\, \ldots \,  j_1 j_2 j_3 \, ; \, a_1 b_1 c_1 d_1 \,  \ldots } \ .
\end{align*}
\end{definition}

The same computation as in Lemma \ref{SymmetriesN} shows that elements in this space fulfil the symmetry:
\begin{equation}\label{SimetriaSimetricos}
 T^{\, \ldots \,  j_1 j_2 j_3 \, ; \, a_1 b_1 c_1 d_1 \,  \ldots } \, = \, T^{\, \ldots \,  c_1 d_1 j_3 \, ; \, a_1 b_1 j_1 j_2 \,  \ldots } \ . 
\end{equation}

\begin{proposition}[\cite{BedetLovelock}]
For any $m \geq 1$: 
$$ S \mathrm{Div}_x^m = 0  \ . $$ 

As a consequence, any natural tensor $T \colon J^2M \to \otimes^{p-3} TX \otimes S^3 TX$ that is divergence-free (in one of the symmetric indices) is indeed a zeroth-order tensor.
\end{proposition}

\proof Due to symmetry (\ref{SimetriaSimetricos}) above:
\begin{align*}
3 \, T^{\, \ldots \,  j_1 j_2 j_3 \, ; \, a_1 b_1 c_1 d_1 \,  \ldots } \, & = \, T^{\, \ldots \,  c_1 d_1 j_3   \, ; \, a_1 b_1 j_1 j_2  \,  \ldots } + T^{\, \ldots \,  c_1 d_1  j_2 \, ; \, a_1 b_1 j_3 j_1 \,  \ldots } + T^{\, \ldots \,  c_1 d_1 j_1 \, ; \, a_1 b_1 j_2 j_3 \,  \ldots } \, = \, 0 \ . 
\end{align*}

\qed


\medskip
If $p = 2k$, with $k \geq 2$, consider the following totally symmetric, natural $2k$-tensors:
$$ S_{2k} := sym ( g^* \otimes \stackrel{k}{\ldots} \otimes g^* ) \ , $$ where $sym$ denotes the symmetrization operator.

As $\nabla g^* = 0$, these tensors $S_{2k}$ are divergence-free.

\begin{theorem}\label{TotSymm} 
If $k \geq 2$, any divergence-free, natural tensor $T \colon J^2M \to S^{2k}TX$ is a constant multiple of $S_{2k}$.

\end{theorem}

\proof Due to the previous Proposition, any totally symmetric, divergence-free tensor  has to be zeroth order. 

By Theorem \ref{Replacement}, the space of zeroth-order natural tensors is isomorphic to $\mathrm{Hom}_{\Or} ( \mathbb{R} , S^{2k} T_xX) = (S^{2k} T_xX)^{\Or}$, so the statement follows from:
$$ \dim_{\mathbb{R}} \left( S^{2k} T_xX \right)^{\Or} = 1 \ , \qquad \forall \, k \in \mathbb{N} \ , $$
which is a trivial computation, using Theorem \ref{InvariantTheory}. 

\qed

\subsubsection*{Non-existence of divergence-free, differential forms}

\begin{definition} For $m , p\geq 1$, let $\Lambda \mathrm{Div}_x^m \subset \Lambda^p T_xX  \otimes  S^m N^2  $
be the vector subspace whose elements satisfy:
\begin{align*}
 0 &= \sum_{(j_p \, c_1 \, d_1 )} T^{\, j_1 \, \ldots \,  j_p \, ; \, a_1 b_1 c_1 d_1 \,  \ldots } \ .
\end{align*}
\end{definition}

\begin{lemma}\label{Formas}
For any $m \geq 1$:
$$ \dim_{\mathbb{R}} \left( \phantom{k^k} \hskip -.4cm \Lambda \mathrm{Div}_x^m \right)^{\Or} = 0  \ . $$
\end{lemma}

\proof  The proof is similar to that of Lemma \ref{CalculoDimension}: let us prove that any total contraction of indices is proportional to:
$$  T^{ii jj \ldots \, ; \, kk \ldots }  $$ and therefore vanishes, because the contraction of two skew-symmetric indices is zero.

We argue by descendent induction on the number of contracted pairs. Given a total contraction, the first index cannot be contracted with the second one (as they are skew-symmetric), but then:
$$ T^{ij \ldots \, ; \, \ldots \, i k \ldots } \, \sim \, T^{kj  \ldots ; \, \ldots \, ii \ldots } $$ and the induction hypothesis applies.

 \qed

\begin{theorem}\label{NoDifferentialForms}
For any $p \geq 1$, there are no divergence-free natural tensors $ \omega \colon J^2 M \to \Lambda^p TX$, but the zero $p$-vector.
\end{theorem}

\proof Due to the previous Lemma, any skew-symmetric tensor that is divergence-free has to be zeroth order. But the space of zeroth-order natural tensors is isomorphic to $$\mathrm{Hom}_{\Or} ( \mathbb{R} , \Lambda^{p} T_xX) = (\Lambda^p T_xX)^{\Or} = 0  \ , $$ because the contraction of two skew-symmetric indices is zero. 

\qed

\begin{remark}
In the oriented case, this result is no longer true: as an example, let $\omega_{P}$ be the Pontryagin 4-form on a six dimensional oriented pseudo-Riemannian manifold $X$.

Its Hodge dual, $* \omega_{P}$, is a natural 2-form that is divergence-free:
$$ \partial * \omega_{P} \sim * ( \d \, \omega_{P} )  \, = \, * ( 0 ) \, = \, 0 \ . $$
\end{remark}

\subsection*{Acknowledgements}

This paper is part of the author's doctoral thesis, written under the supervision of Prof. Juan B. Sancho. 
The author thanks him for his generous advice, as well as A. Gordillo, A. Navarro and J. A. Navarro for helpful comments.


\end{document}